\documentclass[11pt]{article}
\usepackage{latexsym}
                      \usepackage{amssymb}
                      \usepackage{amsmath}
                      \usepackage{amsthm}
                      \usepackage{enumerate}
                      \usepackage{verbatim}
\def\co{\colon\thinspace}

\newtheorem{thm}{Theorem}[section]

\newtheorem{cor}[thm]{Corollary}

\newtheorem{lem}[thm]{Lemma}
\newtheorem{prop}[thm]{Proposition}

\newtheorem{Example}[thm]{Example}
\newenvironment{ex}{\begin{Example}\rm}{\end{Example}}
\newtheorem{Counterexample}[thm]{Counterexample}

\newtheorem{remark}[thm]{Remark}
\newenvironment{rmk}{\begin{remark}\rm}{\end{remark}}
\newtheorem{Fact}[thm]{Fact}

\newtheorem{Nothing}[thm]{$\!\!\!$}

\newcommand{\be}{\begin{equation}}
\newcommand{\ee}{ \end{equation}}
\newcommand{\ba}{\begin{eqnarray}}
\newcommand{\ea}{\end{eqnarray}}
\newcommand{\ban}{\begin{eqnarray*}}
\newcommand{\ean}{\end{eqnarray*}}

\begin{document}
\abovedisplayskip=6pt plus3pt minus3pt \belowdisplayskip=6pt
plus3pt minus3pt
\title{\bf On co-Hopfian nilpotent groups
\rm}
%17B30 Solvable, nilpotent (super) algebras
%20F18 Nilpotent groups
\author{Igor Belegradek
\thanks {Partially supported by NSF Grant \# DMS-0203979.}}
\date{}
\maketitle
\begin{abstract} We characterize co-Hopfian finitely generated 
torsion free nilpotent groups in terms of their Lie algebra
automorphisms, and construct many examples of such groups. 
\end{abstract}
%\tableofcontents

\section{Introduction} 
A group is called {\it co-Hopfian} if it contains no proper
subgroup isomorphic to itself. 
Being co-Hopfian can be interpreted as a rigidity property. 
Many geometrically interesting groups are co-Hopfian, 
including some Kleinian groups~\cite{DP},
some $3$-manifold groups (see~\cite{PW} and references
therein). Any freely indecomposable torsion free 
word-hyperbolic group, other than $\mathbb Z$, is co-Hopfian~\cite{Sel}.
By Mostow rigidity, irreducible lattices in semisimple Lie groups
are co-Hopfian, with the exception of free groups~\cite{Pra}.
The fundamental groups of 
finite volume pinched negatively curved manifolds of
dimension $>2$~\cite{Gro, Bel} are co-Hopfian.
In general, if $i$ is a homotopy invariant of closed manifolds 
which is multiplicative under finite covers
(e.g. Euler characteristic, signature, simplicial volume, or
$L^2$-Betti number), then 
any closed aspherical manifold $M$ with $i(M)\neq 0$ 
has co-Hopfian fundamental group.
 
Following P.~Hall, we call a finitely generated 
torsion free nilpotent group  
an {\it $\mathcal F$-group}. 
An $\mathcal F$-group of class $c$
is called an {\it $\mathcal F_c$-group}.
It seems the only known example of a co-Hopfian $\mathcal F$-group
is due to G.~C.~Smith~\cite[Proposition 4]{Smi} 
(who does not explicitly say the group is co-Hopfian, 
yet the fact is apparent from the proof).

To state the main result we need to review 
some structure theory of $\mathcal F$-groups, mostly due 
to A.~I.~Mal'cev~\cite{Mal2, Mal}. 
Any $\mathcal F$-group $G$ is isomorphic to a cocompact 
discrete subgroup of a simply-connected
nilpotent Lie group $N_G$ whose Lie algebra $L(N_G)$
has rational structure constants. An isomorphism of 
$\mathcal F$-groups $G\to H$ induces an isomorphism
of Lie group $N_G\to N_H$, in particular $N_G$ is determined 
by $G$, up to isomorphism.
A (compact) {\it nilmanifold} is the quotient space of a
simply-connected nilpotent Lie group by a discrete
cocompact subgroup. The fundamental group of a nilmanifold
is an $\mathcal F$-group, and conversely, any 
$\mathcal F$-group is the fundamental group of 
the nilmanifold $N_G/G$.
For nilpotent Lie groups
the exponential map $\exp$ is a diffeomorphism with globally 
defined inverse $\log$. 
The $\mathbb Q$-span of $\log(G)$ is 
a nilpotent Lie subalgebra over $\mathbb Q$, which we denote by $L(G)$.
In fact, $\exp (L(G))$ is a dense subgroup of $N_G$ 
isomorphic to the Mal'cev completion of $G$, and the isomorphism type 
of $L(G)$ depends only on $G$.
Conversely, any nilpotent (finite-dimensional)
Lie algebra $L$ over $\mathbb Q$ arises
as $L(G)$ for some $\mathcal F$-group $G$, and $G$ is 
uniquely determined by $L$ up to commensurability.
This correspondence between $\mathcal F$-groups and nilpotent
Lie algebras over $\mathbb Q$ allows us to translate certain
group-theoretic questions into the Lie algebra language, 
where they are often easier to handle. 
To this end, we obtain the following characterization of
co-Hopfian $\mathcal F$-groups. 

\begin{thm}\label{main} Let $G$ be an $\mathcal F$-group. Then
$G$ is co-Hopfian if and only if $|\det(\Phi)|=1$
for any Lie algebra automorphism $\Phi$ of $L(G)$
that maps $\log(G)$ into itself.
\end{thm} 

In Section~\ref{sec: charnil} we prove, among other things, that 
if $G$ is an $\mathcal F$-group and the Lie algebra $L(G)$ is 
characteristically nilpotent, then $G$ is co-Hopfian.
In Section~\ref{sec: aut} we apply results of 
R.~M.~Bryant and J.~R.~J.~Groves~\cite{BG} to produce
co-Hopfian $\mathcal F$-groups 
with prescribed automorphism groups.
Proposition~\ref{main} is proved in Section~\ref{proofs}, where we
also recall some basic facts on $\mathcal F$-groups.
In Section~\ref{degree} we collect some observations
on maps between nilmanifolds.

\section{Characteristically nilpotent Lie algebras}
\label{sec: charnil}

A Lie algebra $L$ is called {\it characteristically nilpotent}
if each derivation of $L$ is a nilpotent endomorphism.
Any characteristically nilpotent Lie algebra
is nilpotent (see the proof of~\cite[Theorem 1]{LT})
and the class of characteristically nilpotent Lie algebras is closed
under direct sums~\cite{LT}. We get the following.

\begin{cor} \label{charnil}
If $G$ is an $\mathcal F$-group and 
$L(G)$ is characteristically nilpotent, then $G$ is co-Hopfian.
\end{cor}
\begin{proof}
%[Proof of Corollary~\ref{charnil}]
Set $L=L(G)$ and let $\mathrm{Aut}_0(L)$ be the irreducible 
component of the identity in the algebraic group $\mathrm{Aut}(L)$. 
The Lie algebra of $\mathrm{Aut}(L)$ is the algebra of derivations of 
$L$~\cite[II.14, Theorem 16]{Che2},
which is nilpotent since $L$ is characteristically 
nilpotent~\cite[Theorem 1]{LT}.
Hence, each element of $\mathrm{Aut}_0(L)$ is 
unipotent~\cite[V.3.4, Proposition 14]{Che3},
in particular it has determinant $1$. 
Since $\mathrm{Aut}_0(L)$ has finite index in $\mathrm{Aut}(L)$, 
each element of $\mathrm{Aut}(L)$ has determinant $\pm 1$, and 
Proposition~\ref{main} applies. 
\end{proof}

Many explicit examples of characteristically nilpotent Lie algebras
are known (see~\cite{AC} for a survey and references therein).
I am indebted to J.~M.~Ancochea who explained to me that
in each dimension $\ge 7$ there are 
infinitely many pairwise non-isomorphic
characteristically nilpotent Lie algebras over $\mathbb Q$, 
in particular, we deduce
the following.

\begin{cor}
For each $n\ge 7$, the number of
commensurability classes of co-Hopfian $\mathcal F$-groups
of rank $n$ is infinite.
\end{cor}

Here by the {\it rank} of 
an $\mathcal F$-group $G$ we mean $\dim(L(G))$. Incidentally, 
$\dim(L(G))=\dim(N_G)$, and since $N_G$ is contractible and $N_G/G$
is compact, $\dim(N_G)$ is equal to 
the cohomological dimension of $G$. 

Recall that a group $G$ is called
{\it compressible} if any finite index subgroup of $G$
contains a subgroup isomorphic to $G$. 
Some classes of $\mathcal F$-groups
are known to consist of compressible groups~\cite{Smi, Mak}. 
Of course, if an $\mathcal F$-group is compressible, it is not co-Hopfian.
The converse is not true, as the following shows.

\begin{cor}\label{noncompress}
If $L$ is a characteristically nilpotent Lie algebra and
$G$ is an $\mathcal F$-group with $L(G)\cong L$,
then $G\times \mathbb Z$ is not co-Hopfian
and not compressible.
\end{cor}
\begin{proof}
%[Proof of Corollary~\ref{noncompress}]
The group $G\times\mathbb Z$ is not co-Hopfian for it contains
$G\times 2\mathbb Z$. Let $H$ be a proper subgroup of $G$ of
finite index (which exists since $\mathcal F$-groups
are residually finite). 
Arguing by contradiction, assume
$G\times\mathbb Z$ is compressible so that there is an
endomorphism $\phi$ of $G\times\mathbb Z$ with image inside 
$H\times\mathbb Z$, and look at the induced automorphism $\Phi$
of $L(G\times\mathbb Z)\cong L\oplus\mathbb Q$.
Since $L$ is characteristically nilpotent, its center $C(L)$
is contained in $[L, L]$~\cite{LT}.

The derived subalgebra $[L\oplus\mathbb Q, L\oplus\mathbb Q]$
of $L\oplus\mathbb Q$ is just the subspace $[L,L]$.
Since $\Phi$ preserves the derived subalgebra, it
preserves the subspace $C(L)$. 
Since $\exp$ maps the center of the Lie algebra
onto the center of the group, $\phi$ maps $C(G)\times \{0\}$,
where $C(G)$ is the center of $G$, into $C(G)\times\{ 0\}$.

Let $p\co G\times\mathbb Z\to G$
be the projection on the first factor, and
$i\co G\to G\times\mathbb Z$ be the inclusion.
Since $\ker(p)$ lies in the center of $G\times\mathbb Z$
and $\phi$ is injective, $\phi^{-1}(\ker(p))$ also
belongs to the center of $G\times\mathbb Z$,
which is $C(G)\times\mathbb Z$. 
So $i(G)\cap \phi^{-1}(\ker(p))$ is contained in 
$C(G)\times\{ 0\}$, and therefore 
$\phi(i(G))\cap\ker(p)\subseteq C(G)\times\{ 0\}$.
We conclude that $p\circ\phi\circ i$ is an 
injective endomorphism of $G$, whose image lies in $H$,
which is impossible since $G$ is co-Hopfian by 
Corollary~\ref{charnil}.
\end{proof}

Smith noted~\cite{Smi} 
that being compressible is invariant under commensurability
for $\mathcal F$-groups.
It would be interesting to know
whether the same is true for the co-Hopfian property.
In other words, can one decide whether an $\mathcal F$-group
$G$ is co-Hopfian by looking at $L(G)$? 

A related question is whether the co-Hopfian property for $G$
can be read off $L(G)\otimes\mathbb R$, which is the Lie
algebra of $N_G$.
Note that $L(G)$ is characteristically nilpotent if and only if
$L(G)\otimes\mathbb R$ is characteristically nilpotent~\cite{LT},
which immediately implies a stronger version of Corollary~\ref{charnil}.

\begin{cor} 
If $G$ is an $\mathcal F$-group with  $L(G)$ characteristically nilpotent,
then any discrete cocompact subgroup of $N_G$ is co-Hopfian.
\end{cor}

There are several ways to measure the ``size'' of an  
$\mathcal F$-group $G$, the most obvious being the
nilpotency class and the rank. One can 
also measure the ``size'' of $G$ by $rk(G/[G,G])$, which is the first
Betti number of the nilmanifold $N_G/G$. Note that 
$rk(G/[G,G])=\dim(L/[L,L])$ for $L=L(G)$, and $\dim(L/[L,L])$
is equal to the number of 
generators of $L$, as a Lie algebra. In particular, $\dim(L/[L,L])\ge 2$.

Explicit examples of characteristically
nilpotent Lie algebras~\cite{DL, Dye, Fav} 
give rise to the following ``small''
co-Hopfian $\mathcal F$-groups.

\begin{ex}\ \newline
(1) there is a co-Hopfian $\mathcal F_6$-group 
of rank $7$~\cite{DL}. I do not know any 
co-Hopfian $\mathcal F$-groups of rank $<7$. 
Note that nilpotent Lie algebras over 
$\mathbb Q$ of dimension $<7$ are classified in~\cite{Mor}. \newline 
(2) there is a co-Hopfian $\mathcal F_3$-group of
rank $8$~\cite{Fav}. This example is optimal, as far as the class 
is concerned,
since all $\mathcal F$-groups of class $1$ or $2$ 
are compressible~\cite{Smi}. \newline
(3) there is a co-Hopfian $\mathcal F_6$-group $G$ of rank $9$
with $rk(G/[G,G])=2$~\cite{Dye}. This was the first example of
a co-Hopfian $\mathcal F$-group given in~\cite{Smi}.
Note that the construction of~\cite{BG}
discussed in Section~\ref{sec: aut} provides a co-Hopfian 
$\mathcal F$-group of any given value of $rk(G/[G,G])\ge 2$.
\end{ex}

\section{$\mathcal F$-groups with prescribed
automorphism groups}
\label{sec: aut}

A very different source of co-Hopfian $\mathcal F$-groups
is the following construction of Bryant 
and Groves~\cite{BG}. 
Let $n\ge 2$ and $H$ be an arbitrary Zariski-closed subgroup of 
$\mathrm{GL}_n(\mathbb Q)$. 
According to~\cite{BG}, there exists an
$n$-generator nilpotent Lie algebra $L$ over $\mathbb Q$
%(which can be chosen to be solvable of derived length $\le 3$) 
such that $H$ is the image of $\mathrm Aut(L)$ under the homomorphism 
$\mathrm Aut(L)\to\mathrm Aut(L/[L,L])\cong\mathrm{GL}_n(\mathbb Q)$.

For any nilpotent Lie algebra $L$ over $\mathbb Q$,
the kernel of $\mathrm Aut(L)\to\mathrm Aut(L/[L,L])$ consists of 
unipotent automorphisms (see~\cite[pp.136-137]{GS}), 
in particular, each of them
has determinant $1$. 
If in addition $H\subseteq\mathrm{SL}_n(\mathbb Q)$, then
a simple linear algebra argument 
shows that $\det(\phi)=1$ for any $\phi\in\mathrm {Aut}(L)$,
so that by Proposition~\ref{main}, any $\mathcal F$-group $G$
with $L(G)\cong L$ is co-Hopfian. 
(As explained in the proof of Corollary~\ref{arith},
the assumption $H\subseteq\mathrm{SL}_n(\mathbb Q)$ can be dropped
if one is willing to increase $n$ by $1$. Also 
$\mathrm{SL}_n(\mathbb Q)$ can be replaced with the 
group of matrices over $\mathbb Q$ of determinant $\pm 1$).

Recall that a group $A$ is called {\it arithmetic} if there is a
positive integer $n$ and  
$\mathbb Q$-closed subgroup $H$ of $\mathrm{GL}_n(\mathbb C)$
such that $A$ is isomorphic to a subgroup of 
$H\cap\mathrm{GL}_n(\mathbb Q)$ commensurable to 
$H\cap\mathrm{GL}_n(\mathbb Z)$.
It is known that if $G$ is an $\mathcal F$-group $G$, with  
$G/[G,G]$ torsion free, then the image of 
$\mathrm{Aut}(G)$ under the homomorphism 
$\mathrm{Aut}(G)\to\mathrm{Aut}(G/[G,G])$ is arithmetic.
Conversely, it was shown in~\cite{BG} that, 
up to commensurability, any arithmetic subgroup
arises in this way for some $G$. We now note
that in this construction $G$ can be chosen co-Hopfian.

\begin{cor}\label{arith}
For any arithmetic group $A$, there exists a co-Hopfian 
$\mathcal F$-group $G$,
with $G/[G,G]$ torsion free, 
such that the image of $\mathrm{Aut}(G)$
under the homomorphism $\mathrm{Aut}(G)\to\mathrm{Aut}(G/[G,G])$
is commensurable to $A$.
\end{cor}
\begin{proof}
%[Proof of Corollary~\ref{arith}] 
By assumption $A$ is arithmetic in some
$\mathbb Q$-closed subgroup $H$ of $\mathrm{GL}_n(\mathbb C)$.
Since  $\mathrm{GL}_n(\mathbb C)$ can be embedded into 
$\mathrm{SL}_{n+1}(\mathbb C)$ as a $\mathbb Q$-closed subgroup,
we can assume without loss of generality that  
$H$ is a $\mathbb Q$-closed subgroup of 
$\mathrm{SL}_{n+1}(\mathbb C)$.
It remains to repeat the proof of~\cite[Theorem B]{BG},
and note that the resulting group is co-Hopfian by
an argument in the second paragraph of this section.
\end{proof}

\section{Preliminaries and proof of Proposition~\ref{main}}
\label{proofs}

Basic facts about nilpotent groups will be recalled below and 
used without specific reference. The most useful sources for
our purposes 
are~\cite[Chapter 6]{Seg},~\cite[Chapter 5]{CG},~\cite{Mal2, Mal},
and exercises in~\cite{Bou}.

Any $\mathcal F$-group $G$ is embedded in a radicable 
torsion free nilpotent group $G^\mathbb Q$, called the 
{\it Mal'cev completion}
of $G$, that satisfies the following universal property:
any homomorphism of $G$ into a radicable torsion free group $R$
extends to a uniquely determined homomorphism $G^\mathbb Q\to R$.
In particular, any homomorphism of $\mathcal F$-groups $\phi\co G\to H$
extends to a homomorphism of their Mal'cev completions
$\phi^\mathbb Q\co G^\mathbb Q\to H^\mathbb Q$.
Note that if $\phi$ is injective, then so is $\phi^\mathbb Q$,
and the same is true for surjectivity. 
Two $\mathcal F$-groups are commensurable if and only if they have
isomorphic Mal'cev completions. 

The map $\exp\co L(G)\to G^\mathbb Q$
is a bijection with inverse $\log$. 
A map $\psi\co G^\mathbb Q\to H^\mathbb Q$
is a group homomorphism if and only if 
$\log\circ\psi\circ\exp\co L(G)\to L(H)$ is a Lie algebra 
homomorphism.

Let $G$ be an $\mathcal F$-group embedded into $G^\mathbb Q$.
A subgroup $H$ of $G^\mathbb Q$ is called a 
{\it lattice group} if $\log(H)$ is a
subgroup of the vector space $L(G)$ and there exists a set of 
generators of $\log(H)$ that is a basis of the vector space $L(G)$. 
One knows that $G$ is contained in a lattice group
as a subgroup of finite index.
The intersection of all lattice groups containing
$G$ is denoted by $G^\mathrm{lat}$. It follows that $G^\mathrm{lat}$
is a lattice group containing $G$ and $|G^\mathrm{lat}:G|<\infty$. 

\begin{lem} \label{lat-lem} Let $G$ be an $\mathcal F$-group.
If $\phi$ is an injective endomorphism of $G$ and $\phi^\mathbb Q$ 
is the automorphism of $G^\mathbb Q$ extending $\phi$, then \newline
(1) $\phi^\mathbb Q(G^\mathrm{lat})=(\phi(G))^\mathrm{lat}$; \newline
(2) $\phi^\mathbb Q(G^\mathrm{lat})\subseteq G^\mathrm{lat}$ and
$|G^\mathrm{lat}:\phi^\mathbb Q(G^\mathrm{lat})|=|G:\phi(G)|$; \newline
(3) $\phi(G)=G$ if and only if 
$\phi^\mathbb Q(G^\mathrm{lat})=G^\mathrm{lat}$.
\end{lem}
\begin{proof}
(1) and (2) follow from the fact that $\phi^\mathbb Q$ and its inverse
take lattice groups to lattice groups, and (3) follows from (2).
\end{proof}

\begin{cor} Let $G$ be a $\mathcal F$-group.
If $G^\mathrm{lat}$ is co-Hopfian, then so is $G$.
\end{cor}
\begin{proof}
If $\phi$ is an injective endomorphism of $G$,
then $\phi^\mathbb Q$ is an injective endomorphism of $G^\mathrm{lat}$
by Lemma~\ref{lat-lem}(2).
By assumption $\phi^\mathbb Q(G^\mathrm{lat})=G^\mathrm{lat}$,
hence Lemma~\ref{lat-lem}(3) implies that $\phi$ is onto.
\end{proof}

\begin{proof}[Proof of Theorem~\ref{main}]
Let $\phi$ be an injective endomorphism of $G$ and $\phi^\mathbb Q$ 
be the automorphism of $G^\mathbb Q$ extending $\phi$.
Let $\Phi=\log\circ\phi^\mathbb Q\circ\exp$ be the corresponding
automorphism of $L(G)$.
We identify $L(G)$ with $\mathbb Q^n$,
where $n=\dim(L(G))$, in such a way that 
$\log(G^\mathrm{lat})$ corresponds to the standard
$\mathbb Z^n\subset\mathbb Q^n$.
Since $\phi^\mathbb Q$ maps $G^\mathrm{lat}$ into itself, 
in this basis $\Phi$ is represented by a matrix
$F$ that maps $\mathbb Z^n$ into itself, i.e. $F\in M_n(\mathbb Z)$.
Thus, $\phi^\mathbb Q(G^\mathrm{lat})=G^\mathrm{lat}$ if and only if
$F(\mathbb Z^n)=\mathbb Z^n$. Since 
\[F(\mathbb Z^n)=\mathbb Z^n
\Leftrightarrow F\in\mathrm{GL}(\mathbb Z)\Leftrightarrow
|\det(F)|=1\Leftrightarrow |\det(\Phi)|=1,\] 
we use Lemma~\ref{lat-lem}(3) to conclude 
that $\phi(G)=G$ if and only if $|\det(\Phi)|=1$.

Now we are ready to finish the proof. Assume $G$ is co-Hopfian
and let $\Phi$ be an automorphism of $L(G)$ 
that maps $\log(G)$ into itself. Then $\Phi$ defines an 
automorphism of $G^\mathbb Q$ that takes $G$ into itself,
and hence defines an injective endomorphism of $G$
which must be onto since $G$ is co-Hopfian. 
So by above $|\det(\Phi)|=1$.

Conversely, let $\phi$ be an injective endomorphism of $G$.
Then $\phi$ gives rise to a Lie algebra automorphism $\Phi$ of $L(G)$
that maps $\log(G)$ into itself. By assumption $|\det(\Phi)|=1$,
hence by above $\phi$ is onto, and hence $G$ is co-Hopfian.
\end{proof}

\section{Maps of nilmanifolds and generalized Hopfian property}
\label{degree}

An $\mathcal F$-group is co-Hopfian if and only if
the associated nilmanifold $N_G/G$ admits no nontrivial
self-covering map. Understanding self-maps of nilmanifolds
was my original motivation for studying co-Hopfian
$\mathcal F$-groups.

Recall that a group $G$ is called {\it Hopfian} if any epimorphism 
$G\to G$ is an isomorphism. 
It is well known that any $\mathcal F$-group is Hopfian. 
In fact, much more is true: 

\begin{prop}\label{prop: surj}
Any epimorphism of $\mathcal F$-groups of the same rank is
an isomorphism. 
\end{prop}
\begin{proof}
If $\phi\co G\to H$ is a epimorphism of 
$\mathcal F$-groups of the same rank, then 
$\phi$ extends to an epimorphism 
$\phi^\mathbb Q\co  G^\mathbb Q\to H^\mathbb Q$
of the Mal'cev completions. This gives rise
to a Lie algebra epimorphism 
$\log\circ\phi^\mathbb Q\circ\exp\co L(G)\to L(H)$,
which must be injective because $L(G), L(H)$ are finite 
dimensional. It follows that $\phi$ is injective.
\end{proof}

\begin{rmk}
\label{dekimpe}
K.~Dekimpe pointed out to us
that Proposition~\ref{prop: surj} holds for torsion free
virtually polycyclic groups. Indeed, if $\phi\co G\to H$
is an epimorphism of virtually polycyclic groups of the same
(Hirsch) rank, then since the rank is additive under 
extensions, $\ker(\phi)$ has rank $0$, so that $\ker(\phi)$
is finite, and furthermore if $G$ is torsion free, then
$\ker(\phi)$ is trivial. 
\end{rmk}

\begin{prop}
Any nonzero degree map of nilmanifolds is homotopic to a covering map.
In particular, if $N$ is a nilmanifold with co-Hopfian fundamental group,
then any self-map of $N$ of nonzero degree is homotopic to
a diffeomorphism.
\end{prop}
\begin{proof}
Let $f\co M\to N$ be a map of nilmanifolds of nonzero degree $p$.
Let $\tilde f$ be a lift of 
$f$ to the covering $\pi\co\tilde N\to N$ corresponding to the subgroup
$f_*(\pi_1(M))$.
Since degree is multiplicative under composition, $\deg(\pi)$ divides $p$,
in particular $\deg(\pi)$ is finite. Thus 
$f_*\co\pi_1(N)\to f_*(\pi_1(N)))$ is an epimorphism of 
$\mathcal F$-groups of the same rank, so by Proposition~\ref{prop: surj} 
it is an isomorphism.
Hence, $\tilde f$ is homotopic to a diffeomorphism, and so $f$
is homotopic to a covering map.
\end{proof}

\begin{rmk}
M.~Gromov introduced a partial order on the set of closed orientable
manifolds by saying that $M\ge N$ if there is a nonzero degree map
$M\to N$. Thus two nilmanifolds are comparable in this order
if and only if one of the manifolds is a covering space of the other one.
J.~Roitberg~\cite{Roi} constructed two  nilmanifolds
with commensurable fundamental groups that are 
incomparable in this order; the groups are associated with
characteristically nilpotent Lie algebra in~\cite{Dye} and hence
are co-Hopfian.
\end{rmk}

\paragraph{Acknowledgements.}
It is a pleasure to thank J.~M.~Ancochea for showing me 
some examples of characteristically nilpotent Lie algebras,
O.~V.~Belegradek and the referee for comments on the
first version of the paper,
and K.~Dekimpe for Remark~\ref{dekimpe}.

\small
\bibliographystyle{amsplain}
\bibliography{nil}

\

DEPARTMENT OF MATHEMATICS, 253-37, CALIFORNIA INSTITUTE OF TECHNOLOGY,
PASADENA, CA 91125, USA

{\normalsize
{\it email:} \texttt{ibeleg@its.caltech.edu}}
\end{document}